\documentclass[12pt,twoside]{article}
\usepackage[english]{babel}
\usepackage[latin1]{inputenc}
\usepackage{amsmath}
\usepackage{amssymb,amsfonts}
\usepackage{graphicx}                   

\newcommand{\HYP}{\mathbb{H}^3}
\newcommand{\HYN}{\mathbb{H}^n}

\begin{document}
\pagestyle{myheadings}
\markboth{\centerline{Jen\H o Szirmai}}
{The least dense hyperball covering to the regular prism tilings...}
\title
{The least dense hyperball covering to the regular prism tilings \\ in the hyperbolic
$n$-space}

\author{\normalsize{Jen\H o Szirmai} \\
\normalsize Budapest University of Technology and \\
\normalsize Economics Institute of Mathematics, \\
\normalsize Department of Geometry \\
\date{\normalsize{\today}}}

\maketitle


\begin{abstract}
After having investigated the densest packings by congruent hyperballs to the regular prism tilings in the $n$-dimensional hyperbolic space $\HYN$ 
($n \in \mathbb{N},~n \ge 3)$ we consider the dual covering problems and determine the least dense hyperball arrangements and their densities.
\end{abstract}

\newtheorem{theorem}{Theorem}[section]
\newtheorem{corollary}[theorem]{Corollary}
\newtheorem{conjecture}{Conjecture}[section]
\newtheorem{lemma}[theorem]{Lemma}
\newtheorem{exmple}[theorem]{Example}
\newtheorem{defn}[theorem]{Definition}
\newtheorem{rmrk}[theorem]{Remark}
\newenvironment{definition}{\begin{defn}\normalfont}{\end{defn}}
\newenvironment{remark}{\begin{rmrk}\normalfont}{\end{rmrk}}
\newenvironment{example}{\begin{exmple}\normalfont}{\end{exmple}}
\newenvironment{acknowledgement}{Acknowledgement}


\section{Introduction}

In the hyperbolic space $\mathbb{H}^n$ $(n \ge 3)$ a regular prism is the convex hull of two congruent
$(n-1)$ dimensional regular polyhedra in ultraparallel hyperplanes, (i.e. $(n-1)$-planes), related by ,,translation" along the line
joining their centres that is the common perpendicular of the two hyperplanes.
Each vertex of such a tiling is either proper point or every vertex lies on the absolute
quadric of $\mathbb{H}^n$, in this case the prism tiling is called fully asymptotic.
Thus the prism is a polyhedron having at each vertex one $(n-1)$-dimensional regular polytop
and some $(n-1)$-dimensional prisms, meeting at this vertex. \newline
\indent From the definitions of the regular prism tilings and the complete Coxeter orthoschemes of degree $d=1$ (see next Section) follows that
a prism tiling exists in the $n$-dimensional hyperbolic space $\mathbb{H}^n, \ \ (n \geq 3)$ if and only if
exists an appropriate complete Coxeter orthoscheme of degree $d=1$. \newline
\indent The complete Coxeter orthoschemes were classified by {{Im Hof}} in
\cite{IH85} and \cite{IH90} by generalizing the method of {{Coxeter}} and {{B\"ohm}} appropriately.
He showed that they exist only for dimensions $\leq 9$. From this classification follows, that the complete
ortoschemes of degree $d=1$ exist up to 5 dimensions. \newline
\indent The formulas for the hyperbolic covolumes of the considered $n$-dimensional Coxeter tilings are determined in
\cite{K89}, \cite{K91} and \cite{EK}, therefore, it is possible to compute the covolumes of the regular prisms and
the densities of the corresponding ball, horoball and hyperball packings. \newline
\indent In \cite{Sz06-1} and \cite{Sz06-2} we have studied the regular prism tilings and the corresponding optimal hyperball packings in
$\mathbb{H}^n$ $(n=3,4)$ and in the paper  \cite{Sz13} we have extended the in former papers developed method
to 5-dimensional hyperbolic space and construct to each investigated Coxeter tiling a regular prism tiling, 
have studied the corresponding optimal hyperball packings by congruent hyperballs,
moreover, we have determined their metric data and their densities. \newline
\indent In the hyperbolic plane $\mathbb{H}^2$ the universal upper bound of the hypercycle packing density is $\frac{3}{\pi}$
proved by I.~Vermes in \cite{V79} and recently, (to the author's best knowledge) the candidates for the densest hyperball
(hypersphere) packings in the $3,4$ and $5$-dimensional hyperbolic space $\mathbb{H}^n$ are derived by the regular prism
tilings which are studied in papers \cite{Sz06-1}, \cite{Sz06-2} and \cite{Sz13}. \newline
\indent In $\mathbb{H}^2$ the universal lower bound of the hypercycle covering density is $\frac{\sqrt{12}}{\pi}$
determined by I.~Vermes in \cite{V81} but in higher dimensions
there is no result about the hyperball coverings and covering densities. \newline
\indent In this paper we study the $n$-dimensional $(n \ge 3)$ hyperbolic regular prism honeycombs
and the corresponding coverings by congruent hyperballs and we determine their least dense covering densities.
Finally, we formulate a conjecture for the candidate of the least dense hyperball
covering by congruent hyperballs in the 3- and 5-dimensional hyperbolic space.
\section{The projective model and \\ the complete orthoschemes }
We use for $\mathbb{H}^n$ the projective model in the Lorentz space $\mathbb{E}^{1,n}$ of signature $(1,n)$,
i.e.~$\mathbb{E}^{1,n}$ denotes the real vector space $\mathbf{V}^{n+1}$ equipped with the bilinear
form of signature $(1,n)$
$
\langle ~ \mathbf{x},~\mathbf{y} \rangle = -x^0y^0+x^1y^1+ \dots + x^n y^n 
$
where the non-zero vectors
$
\mathbf{x}=(x^0,x^1,\dots,x^n)\in\mathbf{V}^{n+1} \ \  \text{and} \ \ \mathbf{y}=(y^0,y^1,\dots,y^n)\in\mathbf{V}^{n+1},
$
are determined up to real factors, for representing points of $\mathcal{P}^n(\mathbb{R})$. Then $\mathbb{H}^n$ can be interpreted
as the interior of the quadric
$
Q=\{[\mathbf{x}]\in\mathcal{P}^n | \langle ~ \mathbf{x},~\mathbf{x} \rangle =0 \}=:\partial \mathbb{H}^n 
$
in the real projective space $\mathcal{P}^n(\mathbf{V}^{n+1},
\mbox{\boldmath$V$}\!_{n+1})$.

The points of the boundary $\partial \mathbb{H}^n $ in $\mathcal{P}^n$
are called points at infinity of $\mathbb{H}^n $, the points lying outside $\partial \mathbb{H}^n $
are said to be outer points of $\mathbb{H}^n $ relative to $Q$. Let $P([\mathbf{x}]) \in \mathcal{P}^n$, a point
$[\mathbf{y}] \in \mathcal{P}^n$ is said to be conjugate to $[\mathbf{x}]$ relative to $Q$ if
$\langle ~ \mathbf{x},~\mathbf{y} \rangle =0$ holds. The set of all points which are conjugate to $P([\mathbf{x}])$
form a projective (polar) hyperplane
$
pol(P):=\{[\mathbf{y}]\in\mathcal{P}^n | \langle ~ \mathbf{x},~\mathbf{y} \rangle =0 \}.
$
Thus the quadric $Q$ induces a bijection
(linear polarity $\mathbf{V}^{n+1} \rightarrow
\mbox{\boldmath$V$}\!_{n+1})$)
from the points of $\mathcal{P}^n$
onto its hyperplanes.

The point $X [\bold{x}]$ and the hyperplane $\alpha [\mbox{\boldmath$a$}]$
are called incident if $\bold{x}\mbox{\boldmath$a$}=0$ ($\bold{x} \in \bold{V}^{n+1} \setminus \{\mathbf{0}\}, \ \mbox{\boldmath$a$} \in \mbox{\boldmath$V$}_{n+1}
\setminus \{\mbox{\boldmath$0$}\}$).
\begin{definition}
An orthoscheme $\mathcal{S}$ in $\mathbb{H}^n$ $(2\le n \in \mathbb{N})$ is a simplex bounded by $n+1$ hyperplanes $H^0,\dots,H^n$
such that
(see \cite{K91, B--H})
$
H^i \bot H^j, \  \text{for} \ j\ne i-1,i,i+1.
$
\end{definition}

{\it The orthoschemes of degree} $d$ in $\mathbb{H}^n$ are bounded by $n+d+1$ hyperplanes
$H^0,H^1,\dots,H^{n+d}$ such that $H^i \perp H^j$ for $j \ne i-1,~i,~i+1$, where, for $d=2$,
indices are taken modulo $n+3$. For a usual orthoscheme we denote the $(n+1)$-hyperface opposite to the vertex $A_i$
by $H^i$ $(0 \le i \le n)$. An orthoscheme $\mathcal{S}$ has $n$ dihedral angles which
are not right angles. Let $\alpha^{ij}$ denote the dihedral angle of $\mathcal{S}$
between the faces $H^i$ and $H^j$. Then we have
$
\alpha^{ij}=\frac{\pi}{2}, \ \ \text{if} \ \ 0 \le i < j -1 \le n. 
$
The $n$ remaining dihedral angles $\alpha^{i,i+1}, \ (0 \le i \le n-1)$ are called the
essential angles of $\mathcal{S}$.
Geometrically, complete orthoschemes of degree $d$ can be described as follows:
\begin{enumerate}
\item
For $d=0$, they coincide with the class of classical orthoschemes introduced by
{{Schl\"afli}} (see Definitions 2.1).
The initial and final vertices, $A_0$ and $A_n$ of the orthogonal edge-path
$A_iA_{i+1},~ i=0,\dots,n-1$, are called principal vertices of the orthoscheme.
\item
A complete orthoscheme of degree $d=1$ can be interpreted as an
orthoscheme with one outer principal vertex, say $A_n$, which is truncated by
its polar plane $pol(A_n)$ (see Fig.~1 and 3). In this case the orthoscheme is called simply truncated with
outer vertex $A_n$.
\item
A complete orthoscheme of degree $d=2$ can be interpreted as an
orthoscheme with two outer principal vertices, $A_0,~A_n$, which is truncated by
its polar hyperplanes $pol(A_0)$ and $pol(A_n)$. In this case the orthoscheme is called doubly
truncated. We distinguish two different types of orthoschemes but I
will not enter into the details (see \cite{K89}, \cite{K91}).
\end{enumerate}

A $n$-dimensional tiling $\mathcal{P}$ (or solid tessellation, honeycomb) is an infinite set of 
congruent polyhedra (polytopes) fitting together to fill all space $(\mathbb{H}^n~ (n \geqq 2))$ just once, 
so that every face of each polyhedron (polytope) belongs to another polyhedron as well. 
At present the cells are congruent orthoschemes. 
A tiling with orthoschemes exists if and only if the dihedral angle of a tile is a submultiple of $2\pi$ 
(in the hyperbolic plane zero angle is also possible). 

Another approach to describing tilings involves the analysis of their symmetry groups. 
If $\mathcal{P}$ is such a simplex tiling, then any motion taking one cell into another maps the 
whole tiling onto itself. The symmetry group of this tiling is denoted by 
$Sym \mathcal{P}$. 
Therefore the simplex is a fundamental domain of the group $Sym \mathcal{P}$ generated by reflections in its  
$(n-1)$-dimensional hyperfaces.

The scheme of a orthoscheme $S$ is a weighted graph (characterizing $S \subset \mathbb{H}^n$
up to congruence) in which the nodes, numbered by $0,1,\dots,n$ correspond to the bounding hyperplanes of $\mathcal{S}$.
Two nodes are joined by an edge if the corresponding hyperplanes are not orthogonal.
\begin{figure}[ht]
\centering
\includegraphics[width=7cm]{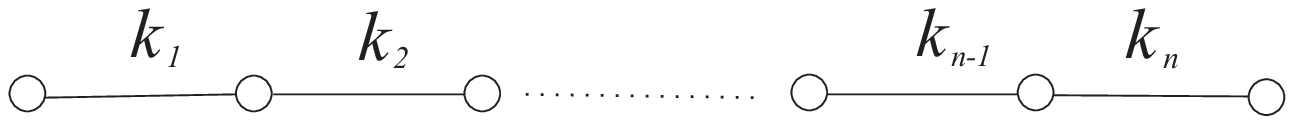}
\end{figure}
For the schemes of {\it complete Coxeter orthoschemes} $\mathcal{S} \subset \mathbb{H}^n$ we adopt the usual conventions and use
them sometimes even in the Coxeter case: If two nodes are related by the weight $\cos{\frac{\pi}{p}}$
then they are joined by a ($p-2$)-fold line for $p=3,~4$ and by a single line marked $p$ for $p \geq 5$.
In the hyperbolic case if two bounding hyperplanes of $S$ are parallel, then the corresponding nodes
are joined by a line marked $\infty$. If they are divergent then their nodes are joined by a dotted line.

The ordered set $[k_1,\dots,k_{n-1},k_n] $ is said to be the 
Coxeter-Schl$\ddot{a}$fli symbol of the simplex tiling $\mathcal{P}$ generated by $\mathcal{S}$. 
To every scheme there is a corresponding
symmetric matrix $(h^{ij})$ of size $(n+1)\times(n+1)$ where $h^{ii}=1$ and, for $i \ne j\in \{0,1,2,\dots,n \}$,
$h^{ij}$ equals $-\cos{\frac{\pi}{k_{ij}}}$ with all angles between the facets $i$,$j$ of $\mathcal{S}$.

For example, $(h^{ij})$ is the so called Coxeter-Schl\"afli matrix of the orthoschem $S$ in the 
5-dimensional hyperbolic space $\mathbb{H}^5$ with
parameters (nodes) $k_1=p,k_2=q,k_3=r,k_4=s,k_5=t$ :
\[
(h^{ij}):=\begin{pmatrix}
1& -\cos{\frac{\pi}{p}}& 0 & 0 & 0 & 0\\
-\cos{\frac{\pi}{p}} & 1 & -\cos{\frac{\pi}{q}}& 0 & 0 &0 \\
0 & -\cos{\frac{\pi}{q}} & 1 & -\cos{\frac{\pi}{r}} & 0 &0 \\
0 & 0 & -\cos{\frac{\pi}{r}} & 1 & -\cos{\frac{\pi}{s}} &0 \\
0 & 0 & 0 & -\cos{\frac{\pi}{s}} & 1&-\cos{\frac{\pi}{t}} \\
0& 0 & 0 & 0 & -\cos{\frac{\pi}{t}} & 1 \\
\end{pmatrix}. \tag{2.4}
\]
\section{Regular prism tilings and their least dense hyperball coverings in $\mathbb{H}^n$}
\subsection{The structure of the $n$-dimensional regular prism tilings}
In hyperbolic space $\mathbb{H}^n$ $(n \ge 3)$ a regular prism is the convex hull of two congruent
$(n-1)$ dimensional regular polyhedra in ultraparallel hyperplanes, (i.e. $(n-1)$-planes), related by
,,translation" along the line joining their centres that is the common perpendicular of the two hyperplanes.
The two regular $(n-1)$-faces of a regular prism are called cover-polytops,
and its other $(n-1)$-dimensional facets are called side-prisms.

In this section we consider the $n$-dimensional regular prism tilings that existence is equivalent to the existence of 
the complete Coxeter orthoschemes of degree $d=1$ that are characterized by their Coxeter-Schl\"afli symbols (see Fig.~1). 
\begin{figure}[ht]
\centering
\includegraphics[width=6.5cm]{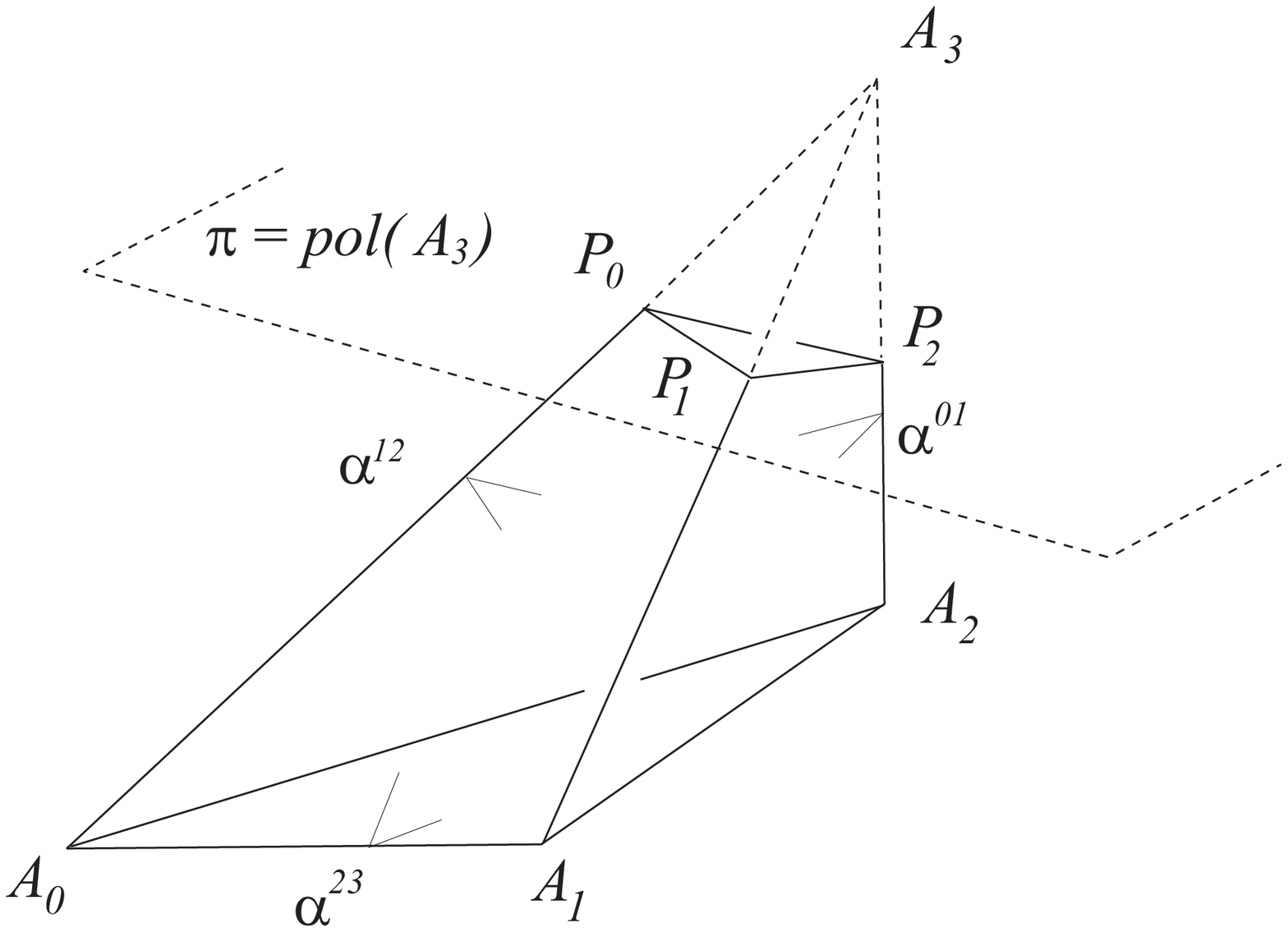} \includegraphics[width=6.5cm]{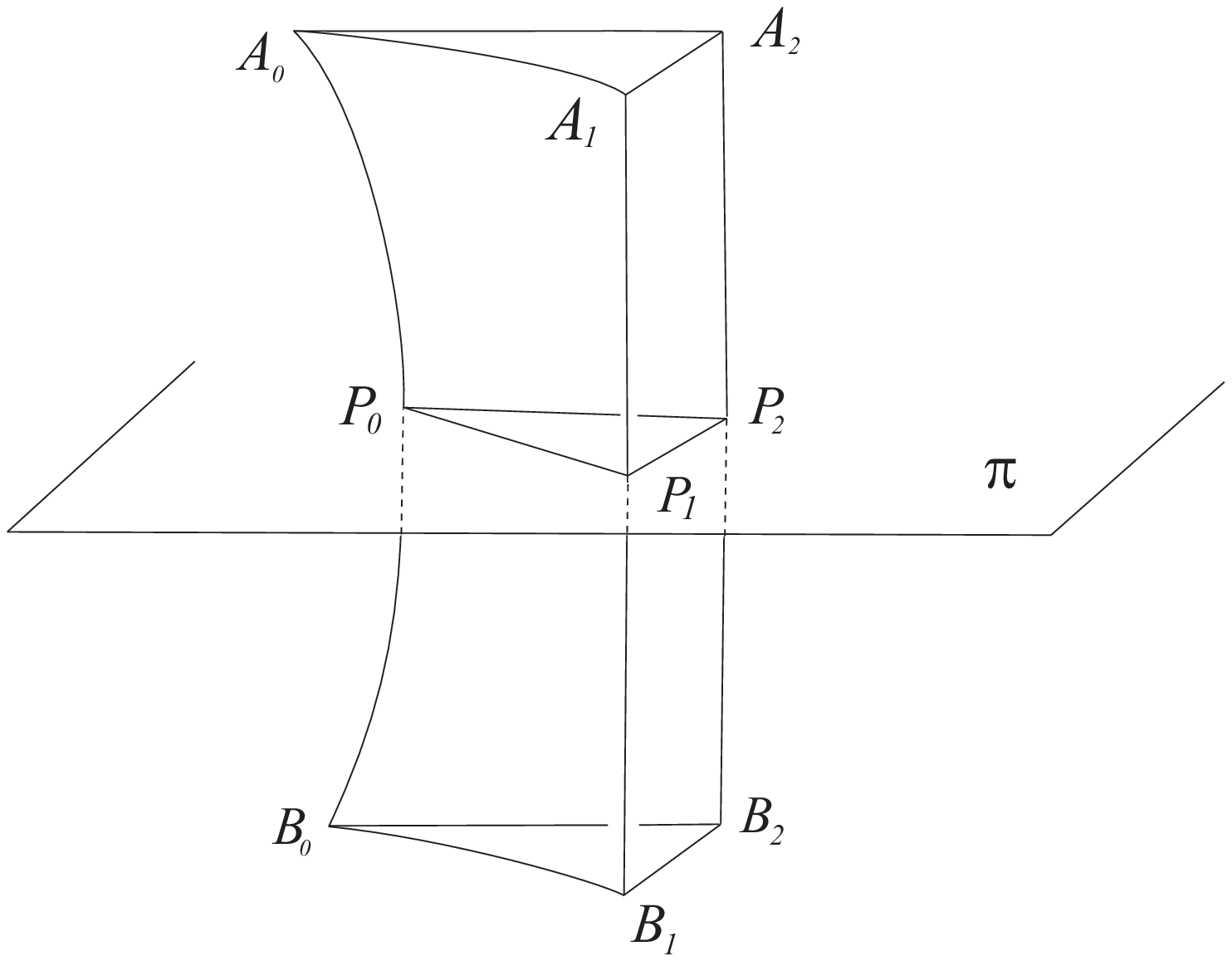}
\caption{A $3$-dimensional complete orthoscheme of degree $d=1$}
\label{}
\end{figure}
The complete Coxeter orthoschemes were classified by {{Im Hof}} in
\cite{IH85} and \cite{IH90} by generalizing the method of {{Coxeter}} and {{B\"ohm}} appropriately.
He showed that they exist only for dimensions $\leq 9$. From this classification follows, that the complete
ortoschemes of degree $d=1$ exist up to 5 dimension. 
\begin{enumerate}
\item In the 3-dimensional hyperbolic space $\HYP$ the regular $p$-gonal non-uniform prism tilings $\mathcal{T}_{pqr}$ with parameters $[p,q,r]$ are the following:
\begin{enumerate}
\item The vertex figure are tetrahedron, octahedron, icosahedron: $ [q,r] = [3,3],~
[3,4],~[3,5]~ \Longrightarrow ~ p>6$,
\item The vertex figure is cube: $ [q,r] = [4,3], ~ \Longrightarrow ~ p>4$,
\item The vertex figure is dodecahedron: $ [q,r] = [5,3], ~ \Longrightarrow ~ p>3$.
\end{enumerate}
\begin{rmrk} 
In the 3-dimensional hyperbolic space there are $3$ infinite series of the totally asymptotic regular prism tilings. 
We do not consider these honeycombs in this paper because hyperball coverings can not be derived from them.
\end{rmrk}
\item The non-uniform compact prism tilings $\mathcal{T}_{pqrs}$ in $\mathbb{H}^4$ with parameters $[p,q,r,s]$ are the following:
\begin{enumerate}
\item $[3,5,3,3]$ : the vertex figure of the tiling is "120-cells": $ [q,r,s] = [5,3,3]$ and the cover faces are
icosahedra $[p,q,r] =[3,5,3].$
\item $[5,3,4,3]$ : the vertex figure of the tiling is "24-cells": $ [q,r,s] = [3,4,3]$ and the cover faces are
dodecahedra $[p,q,r] =[5,3,4].$
\end{enumerate}
\begin{rmrk}
The uniform compact tiling $[4,~3,~3,~5]$ in $\mathbb{H}^4$ is the regular cube honeycomb.
Here the prism is a cube and this tiling is not related with any hyperball packing.
Thus, we do not consider it in this work.
\end{rmrk}
\item The non-uniform compact prism tilings $\mathcal{T}_{pqrst}$ in $\mathbb{H}^5$ with parameters $[p,q,r,s,t]$ are the following:
\begin{enumerate}
\item $[5,3,3,3,3]$ : the vertex figure of the tiling is: $ [q,r,s,t] = [3,3,3,3]$ and the cover faces are
$[p,q,r,s] =[5,3,3,3].$
\item $[5,3,3,4,3]$ : the vertex figure of the tiling is: $ [q,r,s,t] = [3,3,4,3]$ and the cover faces are
$[p,q,r,s] =[5,3,3,4].$
\item $[5,3,3,3,4]$ : the vertex figure of the tiling is: $ [q,r,s,t] = [3,3,3,4]$ and the cover faces are
$[p,q,r,s] =[5,3,3,3].$
\end{enumerate}
\item There is no regular prism tiling in the hyperbolic space $\mathbb{H}^n$, $(n \ge 6)$.
\end{enumerate}
Fig.~3 shows a part of a 5-prism $[p,q,r,s,t]$ where $A_4$ is the centre of a cover-polyhedron,
$A_3$ is the centre of a 3-face of the cover polyhedron, $A_2$ is the midpoint of its $2$-face,
$A_1$ is a midpoint of an edge of this face, and $A_0$ is one vertex (end) of that edge.

Let $B_0,~B_1,~B_2,~B_3,B_4$ be the corresponding points of the other cover-polytop of the regular 5-prism.
The midpoints of the edges $A_iB_i$ $(i\in\{0,1,\dots,4\}$ 
form a hyperplane denoted by $\pi$.
The foot points $P_i ~ (i \in \{ 0,1,2,3,4 \})$
of the perpendiculars dropped from the points
$A_i$ on the plane $\pi$ form the {\it{characteristic (or fundamental) simplex}}
with Coxeter-Schl\"afli symbol $[p,q,r,s]$ (see Fig.~3).

\begin{rmrk}
In $\mathbb{H}^3$ (see \cite{Sz06-1}) the corresponding prisms are called regular $p$-gonal prisms $(p \ge 3)$ in
which the regular polyhedra (the cover-faces)
are regular $p$-gons, and the side-faces are rectangles.
Fig.~1 shows a part of such a prism where $A_2$ is the centre of a regular $p$-gonal
face, $A_1$ is a midpoint of a side of this face, and $A_0$ is one vertex (end) of that side.
Let $B_0,~B_1,~B_2$ be the corresponding points of the other $p$-gonal face of the prism.
\end{rmrk}

Analogously to the $3$-dimensional case, in the $n$-dimensional hyperbolic space $\mathbb{H}^n$ $(n=4,5)$ 
it can be seen that $\mathcal{S}=A_0A_1A_2 \dots A_n$ $P_0P_1P_2 \dots P_n$ is an complete
orthoscheme with degree $d=1$ where $A_n$ is an outer vertex of
$\mathbb{H}^n$ and the points $P_0,P_1,P_2,\dots,P_{n-1}$ lie in its polar hyperplane $\pi$ (see Fig.~3 in the $5$-dimensional hyperbolic space).
The corresponding regular prism $\mathcal{P}$ can be derived by reflections in facets of $\mathcal{S}$ containing the
segment $A_{n-1}P_{n-1}$.
\begin{figure}[ht]
\centering
\includegraphics[width=9cm]{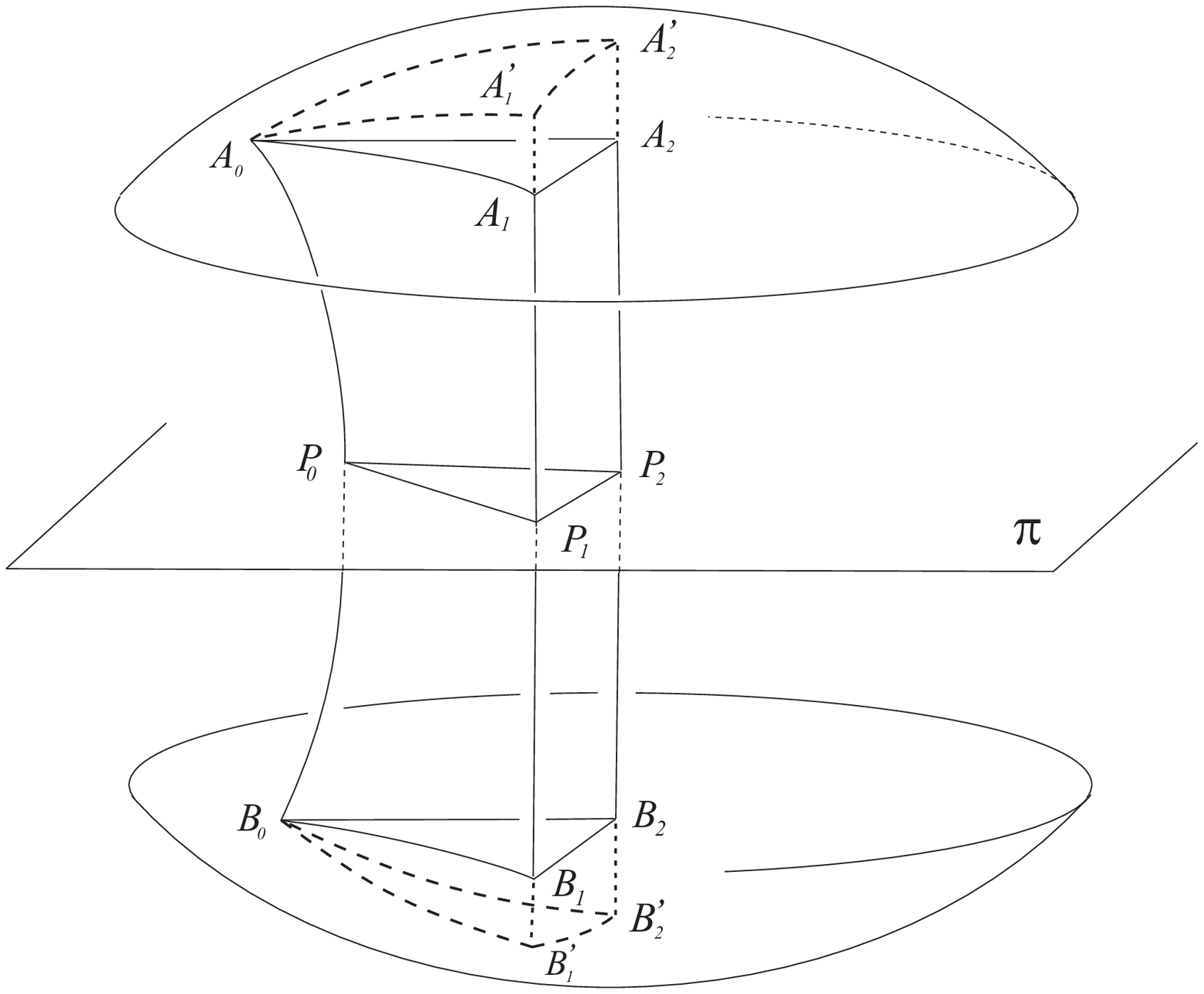}
\caption{}
\label{}
\end{figure}
{\it We consider the images of $\mathcal{P}$ under reflections in its side facets (side prisms).
The union of these $n$-dimensional regular prisms (having the common $\pi$ hyperplane) forms an infinite polyhedron denoted by $\mathcal{F}$.}
$\mathcal{F}$ and its images under reflections in its "cover facets" fill the hyperbolic
space $\mathbb{H}^n$ without overlap and generate $n$-dimensional regular prism tilings $\mathcal{T}$.

From the definitions of the  regular prism tilings and the complete orthoschemes of degree $d=1$ follows that a
regular prism tiling $\mathcal{T}$ exists in the $n$-dimensional hyperbolic space $\mathbb{H}^n, \ \ (n \geq 3)$ if and only if
exists a complete Coxeter orthoscheme of degree $d=1$ with two divergent faces.

On the other hand, if a $n$-dimensional regular prism tiling $[k_1,k_2,\dots k_n]$ exists, then it has to satisfy the following two requirements:
\begin{enumerate}
\item The orthogonal projection of the cover-polytop
on the hyperbolic hyperplane $\pi$ is a regular Coxeter honeycomb with proper vertices and centres.
Using the classical notation of the tesselations, these honeycombs are given by their
Coxeter-Schl\"afli symbols $[k_1,\dots,k_{n-1}]$.
\item The vertex figures about a vertex of
such a prism tiling has to form a $n$-dimensional regular polyhedron.
\end{enumerate}
\subsection{The volumes of the $n$-dimensional \\ Coxeter orthoschemes}

\begin{enumerate}
\item The 3-dimensional hyperbolic space $\HYP$:

{Our polyhedron $A_0A_1A_2P_0P_1P_2$ is a simple frustum orthoscheme with
outer vertex $A_3$ (see Fig.~1) whose volume can be calculated by the following theorem of R.~Kellerhals
(\cite{K89}):}
\begin{theorem}{\rm{(R.~Kellerhals)}} The volume of a three-dimensional hyperbolic
complete ortho\-scheme (except the cases of Lambert cubes) $\mathcal{S}$
is expressed with the essential angles $\alpha_{01},\alpha_{12},\alpha_{23}, \ (0 \le \alpha_{ij} \le \frac{\pi}{2})$
(Fig.~1) in the following form:

\begin{align}
&Vol_3(\mathcal{S})=\frac{1}{4} \{ \mathcal{L}(\alpha_{01}+\theta)-
\mathcal{L}(\alpha_{01}-\theta)+\mathcal{L}(\frac{\pi}{2}+\alpha_{12}-\theta)+ \notag \\
&+\mathcal{L}(\frac{\pi}{2}-\alpha_{12}-\theta)+\mathcal{L}(\alpha_{23}+\theta)-
\mathcal{L}(\alpha_{23}-\theta)+2\mathcal{L}(\frac{\pi}{2}-\theta) \}, \notag
\end{align}
where $\theta \in [0,\frac{\pi}{2})$ is defined by the following formula:
$$
\tan(\theta)=\frac{\sqrt{ \cos^2{\alpha_{12}}-\sin^2{\alpha_{01}} \sin^2{\alpha_{23}
}}} {\cos{\alpha_{01}}\cos{\alpha_{23}}}
$$
and where $\mathcal{L}(x):=-\int\limits_0^x \log \vert {2\sin{t}} \vert dt$ \ denotes the
Lobachevsky function.
\end{theorem}
\begin{figure}[ht]
\centering
\includegraphics[width=8cm]{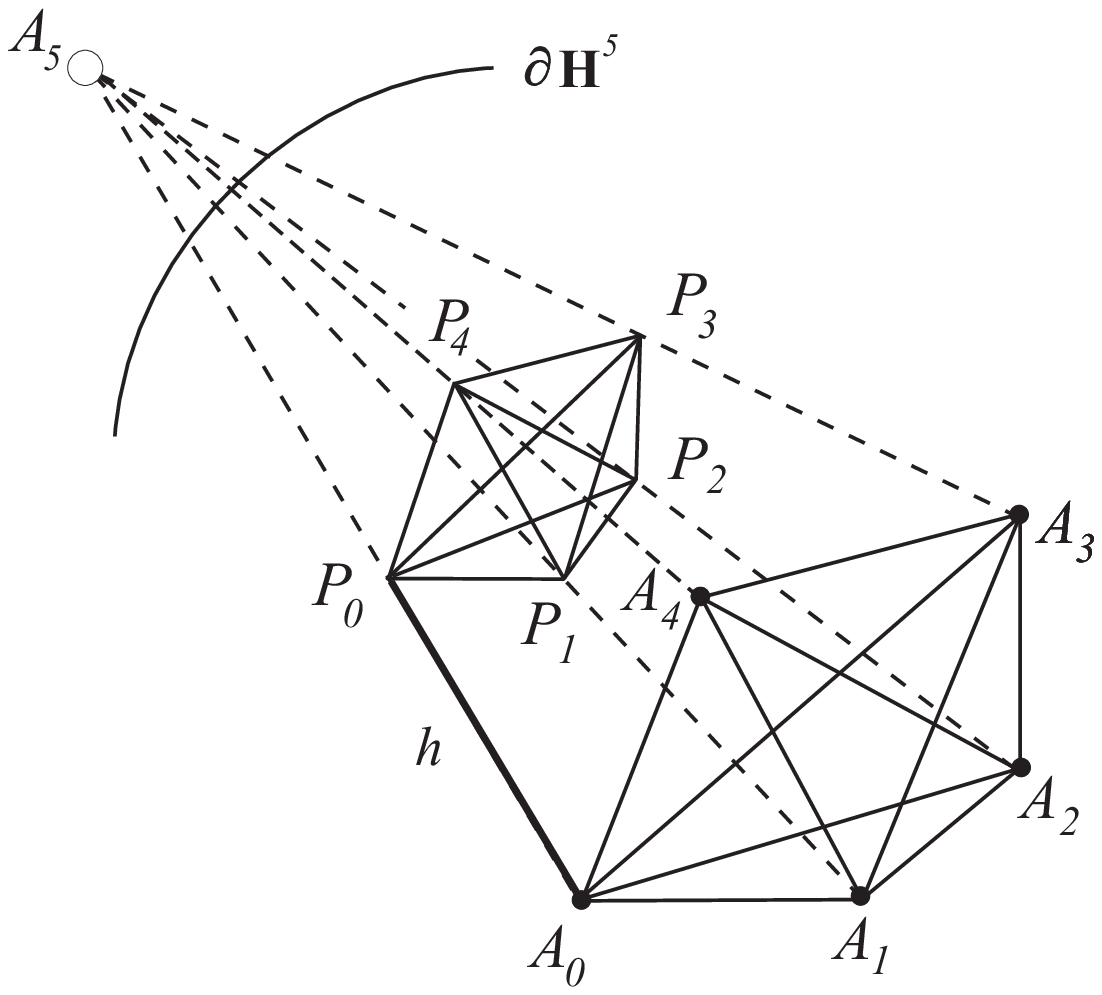}
\caption{}
\label{}
\end{figure}
In our cases for the prism tilings $\mathcal{T}_{pqr}$:~
$\alpha_{01}=\frac{\pi}{p}, \ \ \alpha_{12}=\frac{\pi}{q}, \ \
\alpha_{23}=\frac{\pi}{r}$ .

\item The $4$-dimensional hyperbolic space $\mathbb{H}^4$:

The volumes $Vol_4(S_{pqrs})$ were determined by {{R.~Kellerhals}} in \cite{K91}:
\begin{equation}
Vol_4(S_{3533})=\frac{41 \pi^2}{10800}\approx 0.03746794, \ \ \ Vol_4(S_{5343})=\frac{17 \pi^2}{4320}\approx 0.03883872. \tag{3.1}
\end{equation} 

\item The $5$-dimensional hyperbolic space $\mathbb{H}^5$:

{{R.~Kellerhals}} in \cite{K92} developed a procedure to determine the volumes of 5-dimensionalal hyperbolic orthschemes, 
moreover the volumes of the complete orthoschemes $\mathcal{S}_{pqrst}$ with Coxeter-Schl\"afli symbol $[5,3,3,3,3]$ and 
$[5,3,3,3,4]$ can be computed by volume differential formula of L.~ Schl\"afli with the the following formula (see \cite{EK}):
\begin{equation*}
Vol_5(\mathcal{S}_{pqrst})=\frac{1}{4}\int_{\alpha_i}^{\frac{2\pi}{5}}{Vol_3([5,3,\beta(t)])dt+\frac{\zeta(3)}{3200}} 
\end{equation*}
with a compact tetrahedron $[5,3,\beta(t)]$ whose angle parameter $0<\beta(t)<\frac{\pi}{2}$ is given by
$$
\beta({t})=\arctan{\sqrt{2-\cot^2{t}}}.
$$
Then, the volume of the 3-dimensional orthoscheme face $[5,3,\beta(t)]$ as given by Lobachevsky's formula:
\begin{equation}
\begin{gathered}
Vol_3([5,3,\beta(t)])=\frac{1}{4} \{\mathcal{L}_2\big(\frac{\pi}{5}+\theta(t)\big)-\mathcal{L}_2\big(\frac{\pi}{5}-\theta(t)\big)-\mathcal{L}_2\big(\frac{\pi}{6}+\theta(t)\big)+\\
\mathcal{L}_2\big(\frac{\pi}{6}-\theta(t)\big)+\mathcal{L}_2\big(\beta(t)+\theta(t)\big)-\mathcal{L}_2\big(\beta(t)-\theta(t)\big)+2\mathcal{L}_2\big(\frac{\pi}{2}-\theta(t)\big) \end{gathered} \tag{3.2}
\end{equation}
where $\mathcal{L}(\omega)$ is the Lobachevsky's function,
$
\theta(t)=\arctan\frac{\sqrt{1-4\sin^2\frac{\pi}{5} \sin^2{\beta(t)}}}{2\cos\frac{\pi}{5} \cos\beta(t)}
$ and $\beta(t)=\frac{\pi}{3}$ or $\frac{\pi}{4}$.
\end{enumerate}
\subsection{The least dense hyperball coverings}
The equidistance surface (or hypersphere) is a quadratic surface at a constant distance
from a plane in both halfspaces. The infinite body of the hypersphere is called hyperball.

The $n$-dimensional {\it half hypersphere } $(n=3,4,5)$ with distance $h$ to a hyperplane $\pi$ is denoted by $\mathcal{H}_n^h$. 
The volume of a bounded hyperball piece $\mathcal{H}_n^h(\mathcal{A}_{n-1})$ 
delimited by a $(n-1)$-polytop $\mathcal{A}_{n-1} \subset \pi$, $\mathcal{H}_n^h$ and by some to $\pi$ orthogonal 
hyperplanes derived by the facets of $\mathcal{A}_{n-1}$ can be determined by the formulas (3.3), (3,4) and (3.5) that follow by the generalization of the 
classical method of {{J.~Bolyai}}:
\begin{equation}
Vol_3(\mathcal{H}_3^h(\mathcal{A}_2))=\frac{1}{4}Vol_2(\mathcal{A}_{2})\left[k \sinh \frac{2h}{k}+ 
2 h \right], \tag{3.3}
\end{equation}
\begin{equation}
Vol_4(\mathcal{H}_4^h(\mathcal{A}_3))=\frac{1}{8} Vol_3(\mathcal{A}_{3})k \left[ \frac{2}{3} \sinh \frac{3h}{k}+
6 \sinh \frac{h}{k} \right], \tag{3.4}
\end{equation}
\begin{equation}
Vol_5(\mathcal{H}_5^h(\mathcal{A}_4))=\frac{1}{16} Vol_4(\mathcal{A}_4)\left[ k  \left( \frac{1}{2} \sinh \frac{4h}{k}+
4 \sinh \frac{2h}{k}\right) +6 h \right], \tag{3.5}
\end{equation}
where the volume of the hyperbolic $(n-1)$-polytop $\mathcal{A}_{n-1}$ lying in the plane $\pi$ is $Vol_{n-1}(\mathcal{A}_{n-1})$.
The constant $k =\sqrt{\frac{-1}{K}}$ is the natural length unit in
$\mathbb{H}^n$. $K$ will be the constant negative sectional curvature.

We consider one from the former described $n$-dimensional $(n=3,4,5)$ regular prism tilings $\mathcal{T}$ and
the corresponding infinite polyhedron $\mathcal{F}$ derived from that (the union of $n$-dimensional regular prisms having 
the common hyperplane $\pi$). $\mathcal{F}$ and its images under reflections in its "cover facets" fill the hyperbolic
space $\mathbb{H}^n$ without overlap. 

If we start with a given congruent hyperball covering in the $n$-dimensional hyperbolic space $\mathbb{H}^n$ and shrink the heights of hyperballs
until they finally do not cover the space any more, then the 
minimal height (radius) defines the least dense covering to a given hyperball arrangement.
The thresfold value is called {\it the minimal covering height (radius)} of the given hyperball arrangement.

{\it{We are looking for the smallest half hyperball ${\mathcal{H}_n^{h}}$ containing $\mathcal{F}$ with minimal covering height.}}

The smallest half hypersphere $\mathcal{H}_n^h$ contains the cover-faces of the regular $n$-prisms containing by $\mathcal{F}$.
Therefore, the minimal distance from the $(n-1)$-midplane $\pi$ will be $h=P_0A_0>P_iA_i$ $(i\in\{1,2,3,4\}$ (Fig.~1 and Fig.~3).

The smallest hypersphere $\mathcal{H}_n^h$ covers the infinite polyhedron $\mathcal{F}$ with minimal covering height 
thus, we obtain by the images of $\mathcal{H}_{n}^{h}$  locally least dense hyperball covering to the tiling $\mathcal{T}$.
\subsection{The minimal covering height of the hyperballs to regular prism tilings}
The points $P_0[{\mathbf{p}}_0]$ and $A_0[{\mathbf{a}}_0]$ are proper points of the hyperbolic $n$-space and
$P_0$ lies on the polar hyperplane $pol(A_n)[\mbox{\boldmath$a$}^n]$ of the outer point $A_n$ thus
\begin{equation}
\begin{gathered}
\mathbf{p}_0 \sim c \cdot \mathbf{a}_n+\mathbf{a}_0 \in \mbox{\boldmath$a$}^n \Leftrightarrow
c \cdot \mathbf{a}_n \mbox{\boldmath$a$}^n+\mathbf{a}_0 \mbox{\boldmath$a$}^n=0 \Leftrightarrow
c=-\frac{\mathbf{a}_0 \mbox{\boldmath$a$}^n}{\mathbf{a}_n \mbox{\boldmath$a$}^n} \Leftrightarrow \\
\mathbf{p}_0 \sim -\frac{\mathbf{a}_0 \mbox{\boldmath$a$}^n}{\mathbf{a}_n \mbox{\boldmath$a$}^n}
\mathbf{a}_n+\mathbf{a}_0 \sim \mathbf{a}_0 (\mathbf{a}_n \mbox{\boldmath$a$}^n) - \mathbf{a}_n (\mathbf{a}_0 \mbox{\boldmath$a$}^n)=
\mathbf{a}_0 h_{nn}-\mathbf{a}_n h_{0n},
\end{gathered} \tag{3.6}
\end{equation}
where $h_{ij}$ is the inverse of the Coxeter-Schl\"afli matrix $c^{ij}$ (see (2.4)) of the orthoscheme $\mathcal{S}$.
The hyperbolic distance $h$ can be calculated by the following formula \cite{M89}:
\[
\begin{gathered}
\cosh{P_0A_0}=\cosh{h}=\frac{- \langle {\mathbf{p}}_0, {\mathbf{a}}_0 \rangle }
{\sqrt{\langle {\mathbf{p}}_0, {\mathbf{p}}_0 \rangle \langle {\mathbf{a}}_0, {\mathbf{a}}_0 \rangle}}= \\ =\frac{h_{0n}^2-h_{00}h_{nn}}
{\sqrt{h_{00}\langle \mathbf{p}_0, \mathbf{p}_0 \rangle}} =
\sqrt{\frac{h_{00}~h_{nn}-h_{0n}^2}
{h_{00}~h_{nn}}}.
\end{gathered} \tag{3.7}
\]
The volume of the polyhedron (complete orthoscheme of degree 1) $\mathcal{S}$ is denoted by $Vol_n(\mathcal{S})$.

For the density of the covering it is sufficient to relate the volume of the minimal covering hyperball piece to that of
corresponding polyhedron $\mathcal{S}$ (see Fig.~2 and 3) because the tiling can be constructed of such polyhedron. This polytope
and its images in $\mathcal{F}$ divide the $\mathcal{H}_n^{h}$ into congruent horoball
pieces whose volume is denoted by $Vol_n({\mathcal{H}_{n}^{h}(\mathcal{A}_{n-1})})$. We illustrate in the 3-dimensional case such a
hyperball piece $A_0 A_1' A_2' P_0 P_1 P_2$ in Fig.~2.

The density of the least dense hyperball covering to the $n$-dimensional regular prism tiling $\mathcal{T}$ $(n=3,4,5)$ 
is defined by the following formula:
\begin{definition}
\begin{equation}
\delta^{min}(\mathcal{T}):=\frac{Vol_n(\mathcal{H}_n^h(\mathcal{A}_{n-1}))}{Vol_n({\mathcal{S}})}. \tag{3.8}
\end{equation}
\end{definition}
\pagebreak
\section{The data of the hyperball coverings}
\begin{enumerate}
\item 3-dimensional hyperbolic space $\mathbb{H}^3$

By the formulas (3.3), (3.7), (3.8) and by the Theorem 3.4 we can calculate the data
and the densities of the least dense hyperball coverings to each
regular prism tiling in the hyperbolic space $\mathbb{H}^3$ which are summarized in Tables 1-5.

For every prism tiling we have determined in a suitable interval
the graph of the functions $Vol_3(\mathcal{H}_3^{h})(p)$ and
$\delta^{min}_{pqr}(p)$ as continuous functions of $p$ with fixed $q,~r$.
In Fig.~4 we have described these functions for the case $[p,~3,~3]$.

From the formulas (3.7) follows that the function $h(p)$ is decreasing and the function $Vol_3(\mathcal{H}_3^h)$ is increasing
in its domain of definition. By the Theorem 3.4 and the formula (3.3) it can be seen that the function $Vol_3(\mathcal{S}_{pqr})(p)$ increases 
similarly to the function
$\delta^{min}_{pqr}(p)$ in their domains of definition.
\begin{figure}[ht]
\centering
\includegraphics[width=6.5cm]{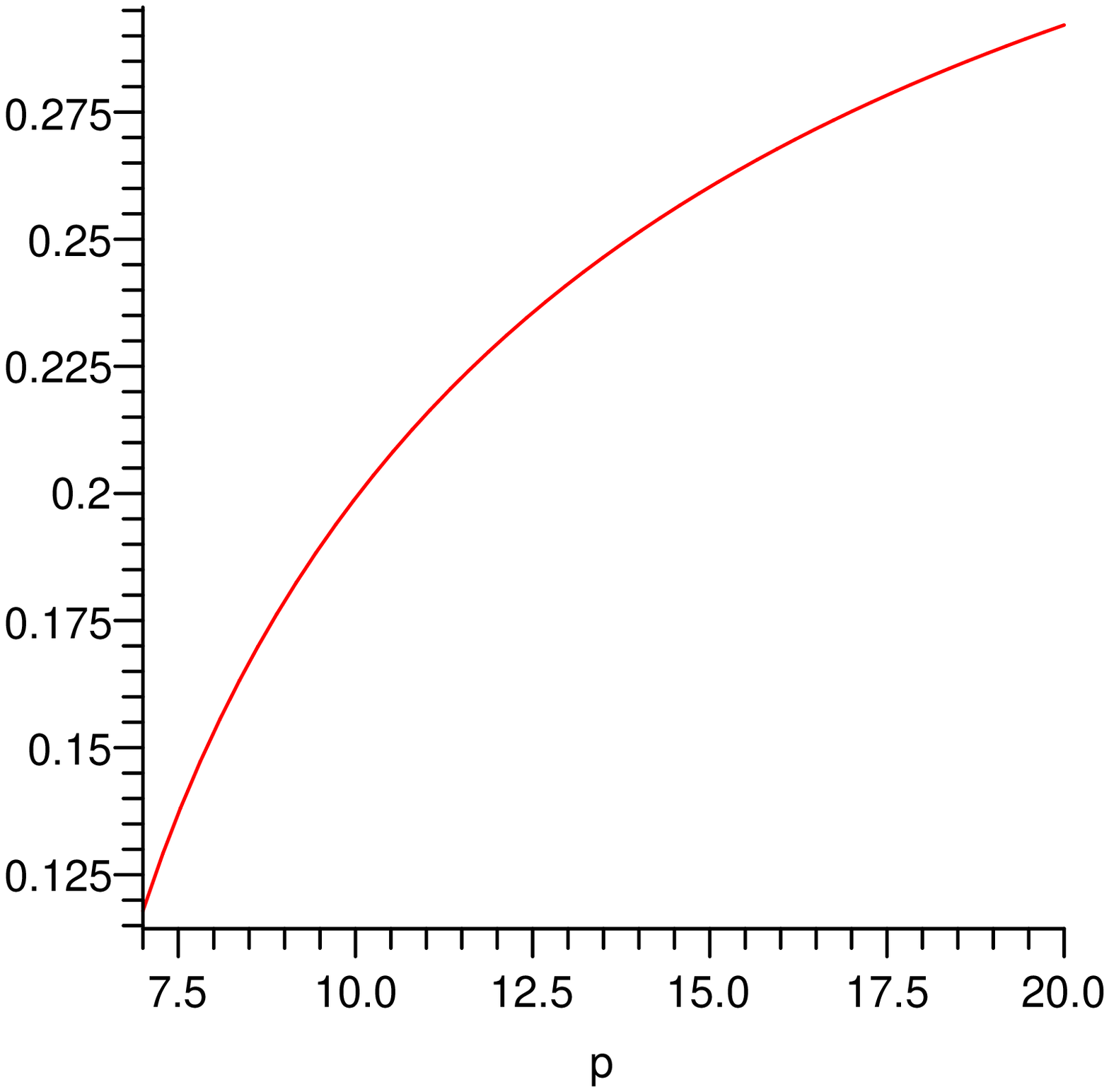} \includegraphics[width=6.5cm]{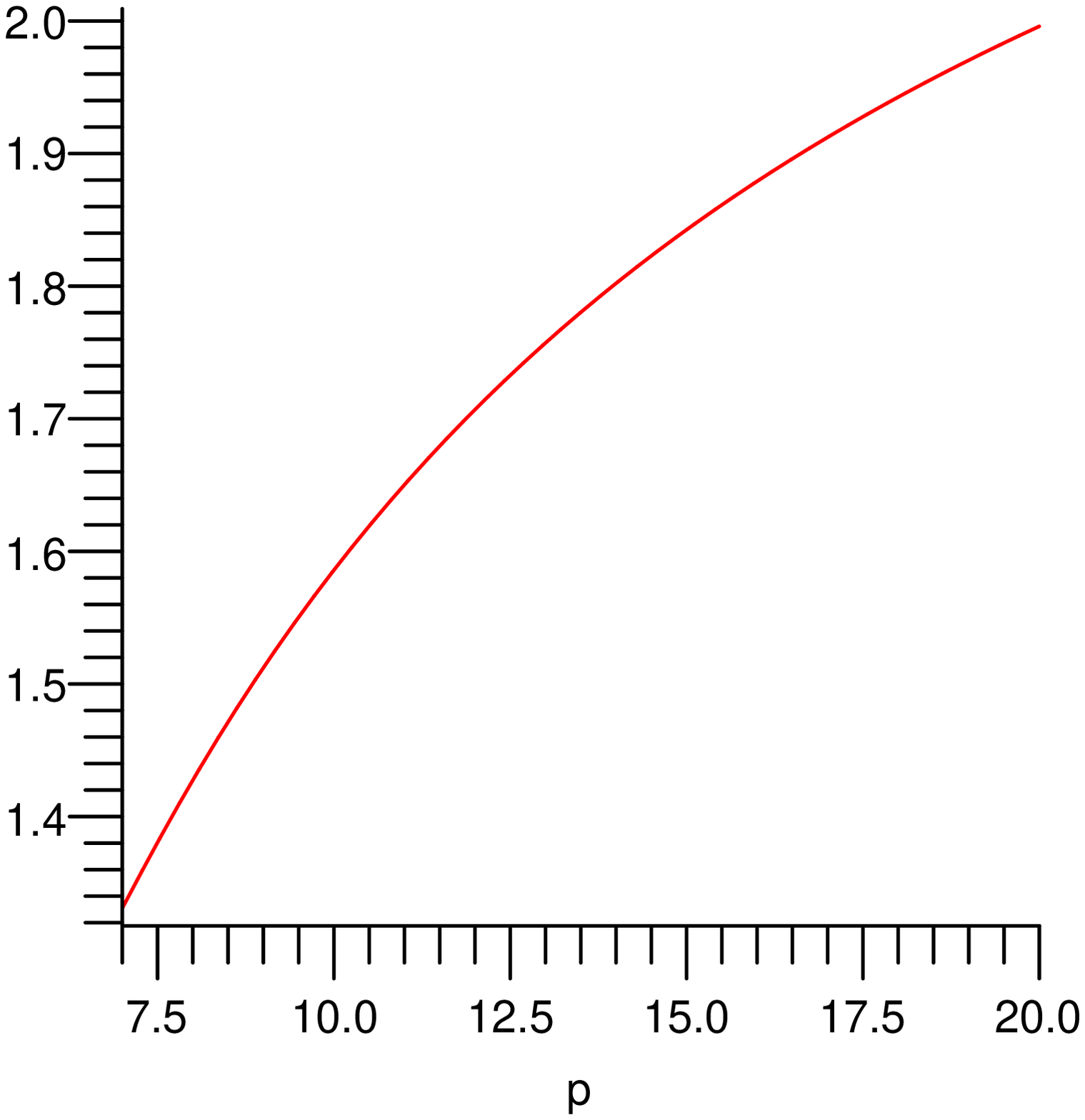}
\caption{The functions $Vol_3(\mathcal{H}_3^h)(p)$ and $\delta_{pqr}^{min}(p)$ to parameters $q=3,~r=3$}
\label{}
\end{figure}
We have determined in \cite{Sz06-1} for each case the limit of the function $Vol_3(\mathcal{S}_{pqr})(p)$ as
$p \to \infty$ (see Table 1-5).
\medbreak
\footnotesize
\centerline{\vbox{
\halign{\strut\vrule~\hfil $#$ \hfil~\vrule
&\quad \hfil $#$ \hfil~\vrule
&\quad \hfil $#$ \hfil\quad\vrule
&\quad \hfil $#$ \hfil\quad\vrule
&\quad \hfil $#$ \hfil\quad\vrule
\cr
\noalign{\hrule}
\multispan5{\strut\vrule\hfill\bf Table 1, $q=3$, $r=3$  \hfill\vrule}%
\cr
\noalign{\hrule}
\noalign{\vskip2pt}
\noalign{\hrule}
p & h & Vol_3(\mathcal{S}_{pqr}) & Vol_3(\mathcal{H}_3^h) & \delta_{pqr}^{min} \cr
\noalign{\hrule}
7 & 1.06738516 & 0.08856157 & 0.11786931 & 1.33093073 \cr
\noalign{\hrule}
8 & 0.89197684 & 0.10721273 & 0.15304272 & 1.42746787 \cr
\noalign{\hrule}
9 & 0.81695936 & 0.11824897 & 0.17882183 & 1.51224843 \cr
\noalign{\hrule}
\vdots & \vdots  & \vdots  & \vdots  & \vdots \cr
\noalign{\hrule}
20 & 0.68135915 & 0.14636009 & 0.29212819 & 1.99595522 \cr
\noalign{\hrule}
\vdots & \vdots  & \vdots  & \vdots  & \vdots \cr
\noalign{\hrule}
50 & 0.66193428 & 0.15167070 & 0.35361352 & 2.33145566 \cr
\noalign{\hrule}
\vdots & \vdots  & \vdots  & \vdots  & \vdots \cr
\noalign{\hrule}
100 & 0.65933594 & 0.15241399 & 0.37580045 & 2.46565592 \cr
\noalign{\hrule}
p \to \infty & 0.65847895 & 0.15266093 & 0.39911431 & 2.61438404 \cr
\noalign{\hrule}}}}
\smallbreak
\centerline{\vbox{
\halign{\strut\vrule~\hfil $#$ \hfil~\vrule
&\quad \hfil $#$ \hfil~\vrule
&\quad \hfil $#$ \hfil\quad\vrule
&\quad \hfil $#$ \hfil\quad\vrule
&\quad \hfil $#$ \hfil\quad\vrule
\cr
\noalign{\hrule}
\multispan5{\strut\vrule\hfill\bf Table 2, $q=4$, $r=3$  \hfill\vrule}%
\cr
\noalign{\hrule}
\noalign{\vskip2pt}
\noalign{\hrule}
p & h & Vol_3(\mathcal{S}_{pqr}) & Vol_3(\mathcal{H}_{3}^h) & \delta_{pqr}^{min} \cr
\noalign{\hrule}
5 & 1.16973604 & 0.16596371 & 0.29370599 & 1.76970002 \cr
\noalign{\hrule}
6 & 1.03171853 & 0.19616337 & 0.38853781 & 1.98068483 \cr
\noalign{\hrule}
7 & 0.97766375 & 0.21217704 & 0.45589711 & 2.14866375 \cr
\noalign{\hrule}
\vdots & \vdots  & \vdots  & \vdots  & \vdots \cr
\noalign{\hrule}
20 & 0.89038430 & 0.24655736 & 0.73257578 & 2.97121844 \cr
\noalign{\hrule}
\vdots & \vdots  & \vdots  & \vdots  & \vdots \cr
\noalign{\hrule}
50 & 0.88277651 & 0.25026133 & 0.83138640 & 3.32207298 \cr
\noalign{\hrule}
\vdots & \vdots  & \vdots  & \vdots  & \vdots \cr
\noalign{\hrule}
100 & 0.88172298 & 0.25078571 & 0.86594296 & 3.45291981 \cr
\noalign{\hrule}
p \to \infty & 0.88137359 & 0.25096025 & 0.901474965 & 3.59210258 \cr
\noalign{\hrule}}}}
\smallbreak
\centerline{\vbox{
\halign{\strut\vrule~\hfil $#$ \hfil~\vrule
&\quad \hfil $#$ \hfil~\vrule
&\quad \hfil $#$ \hfil\quad\vrule
&\quad \hfil $#$ \hfil\quad\vrule
&\quad \hfil $#$ \hfil\quad\vrule
\cr
\noalign{\hrule}
\multispan5{\strut\vrule\hfill\bf Table 3, $q=3$, $r=4$  \hfill\vrule}%
\cr
\noalign{\hrule}
\noalign{\vskip2pt}
\noalign{\hrule}
p & h & Vol_3(\mathcal{S}_{pqr}) & Vol_3(\mathcal{H}_{3}^h) & \delta_{pqr}^{min} \cr
\noalign{\hrule}
7 & 1.67069036 & 0.16297337 & 0.32636556 & 2.00256989 \cr
\noalign{\hrule}
8 & 1.45701769 & 0.18789693 & 0.39605076 & 2.10780855 \cr
\noalign{\hrule}
9 & 1.36099521 & 0.20295023 & 0.44917804 & 2.21324236 \cr
\noalign{\hrule}
\vdots & \vdots  & \vdots  & \vdots  & \vdots \cr
\noalign{\hrule}
20 & 1.17836243 & 0.24206876 & 0.69524288 & 2.87208836 \cr
\noalign{\hrule}
\vdots & \vdots  & \vdots  & \vdots  & \vdots \cr
\noalign{\hrule}
50 & 1.15109755 & 0.24956032 & 0.83516608 & 3.34654595 \cr
\noalign{\hrule}
\vdots & \vdots  & \vdots  & \vdots  & \vdots \cr
\noalign{\hrule}
100 & 1.14742750 & 0.25061105 & 0.88666316 & 3.53800507 \cr
\noalign{\hrule}
p \to \infty & 1.14621583 & 0.25096025 & 0.94135352 & 3.75100645 \cr
\noalign{\hrule}}}}
\smallbreak
\centerline{\vbox{
\halign{\strut\vrule~\hfil $#$ \hfil~\vrule
&\quad \hfil $#$ \hfil~\vrule
&\quad \hfil $#$ \hfil\quad\vrule
&\quad \hfil $#$ \hfil\quad\vrule
&\quad \hfil $#$ \hfil\quad\vrule
\cr
\noalign{\hrule}
\multispan5{\strut\vrule\hfill\bf Table 4, $q=5$, $r=3$  \hfill\vrule}%
\cr
\noalign{\hrule}
\noalign{\vskip2pt}
\noalign{\hrule}
p & h & Vol_3(\mathcal{S}_{pqr}) & Vol_3(\mathcal{H}_{3}^h) & \delta_{pqr}^{min} \cr
\noalign{\hrule}
4 & 1.59191259 & 0.21298841 & 0.59818156 & 2.80851695 \cr
\noalign{\hrule}
5 & 1.40035889 & 0.26319948 & 0.86382502 & 3.28201651 \cr
\noalign{\hrule}
6 & 1.34187525 & 0.28635531 & 1.04400841 & 3.64584964 \cr
\noalign{\hrule}
\vdots & \vdots  & \vdots  & \vdots  & \vdots \cr
\noalign{\hrule}
20 & 1.26417766 & 0.32848945 & 1.71902032 & 5.23310655 \cr
\noalign{\hrule}
\vdots & \vdots  & \vdots  & \vdots  & \vdots \cr
\noalign{\hrule}
50 & 1.25939250 & 0.33171659 & 1.90998817 & 5.75789160 \cr
\noalign{\hrule}
\vdots & \vdots  & \vdots  & \vdots  & \vdots \cr
\noalign{\hrule}
100 & 1.25872455 & 0.33217467 & 1.97599824 & 5.94867215 \cr
\noalign{\hrule}
p \to \infty & 1.25850276 & 0.33232721 & 2.04337965 & 6.14869802 \cr
\noalign{\hrule}}}}
\smallbreak
\centerline{\vbox{
\halign{\strut\vrule~\hfil $#$ \hfil~\vrule
&\quad \hfil $#$ \hfil~\vrule
&\quad \hfil $#$ \hfil\quad\vrule
&\quad \hfil $#$ \hfil\quad\vrule
&\quad \hfil $#$ \hfil\quad\vrule
\cr
\noalign{\hrule}
\multispan5{\strut\vrule\hfill\bf Table 5, $q=3$, $r=5$  \hfill\vrule}%
\cr
\noalign{\hrule}
\noalign{\vskip2pt}
\noalign{\hrule}
p & h & Vol_3(\mathcal{S}_{pqr}) & Vol_3(\mathcal{H}_{3}^h) & \delta_{pqr}^{min} \cr
\noalign{\hrule}
7 & 2.26142836 & 0.23325784 & 0.94559299 & 4.05385304 \cr
\noalign{\hrule}
8 & 2.03433214 & 0.26094396 & 1.08972405 & 4.17608465 \cr
\noalign{\hrule}
9 & 1.93012831 & 0.27782716 & 1.20377506 & 4.33291999 \cr
\noalign{\hrule}
\vdots & \vdots  & \vdots  & \vdots  & \vdots \cr
\noalign{\hrule}
20 & 1.726831092 & 0.32216770 & 1.76349123 & 5.47383009 \cr
\noalign{\hrule}
\vdots & \vdots  & \vdots  & \vdots  & \vdots \cr
\noalign{\hrule}
50 & 1.69577933 & 0.33072584 & 2.10004252 & 6.34979864 \cr
\noalign{\hrule}
\vdots & \vdots  & \vdots  & \vdots  & \vdots \cr
\noalign{\hrule}
100 & 1.69158357 & 0.33192770 & 2.22690186 & 6.70899675 \cr
\noalign{\hrule}
p \to \infty & 1.69019748 & 0.33232721 & 2.36333702 & 7.11147614 \cr
\noalign{\hrule}}}}
\normalsize
\item 4-dimensional hyperbolic space $\mathbb{H}^4$

We obtain the densities of the least dense hyperball packing to the regular prism tilings by the results (3.1) and by the formulas (3.4)and (3.7)
which are summarized in the Table 6.
\footnotesize
\medbreak
\centerline{\vbox{
\halign{\strut\vrule~\hfil $#$ \hfil~\vrule
&\quad \hfil $#$ \hfil~\vrule
&\quad \hfil $#$ \hfil\quad\vrule
&\quad \hfil $#$ \hfil\quad\vrule
&\quad \hfil $#$ \hfil\quad\vrule
\cr
\noalign{\hrule}
\multispan5{\strut\vrule\hfill\bf Table 6, 4-dimensional cases \hfill\vrule}%
\cr
\noalign{\hrule}
\noalign{\vskip2pt}
\noalign{\hrule}
 \mathcal{T} & h & Vol_4(\mathcal{S}_{pqrs}) & Vol_4(\mathcal{H}_{4}^{h}) & \delta_{pqrs}^{min} \cr
\noalign{\hrule}
[3,5,3,3] & 1.96333162 & \frac{41 \pi^2}{10800} & 0.69028590 & 18.42337348 \cr
\noalign{\hrule}
[5,3,4,3] & 1.46935174 & \frac{17 \pi^2}{4320} & 0.178146199 & 4.58681940 \cr
\noalign{\hrule}
\noalign{\hrule}}}}
\normalsize
\item 5-dimensional hyperbolic space $\mathbb{H}^5$

The date of the least dense hyperball packing to the regular prism tilings can be detrmined by the formulas (3.2), (3.5) and (3.7) 
which are summarized in the Table 7.
\begin{rmrk} 
In the 5-dimensional hyperbolic space there is a totally asymptotic regular prism tiling $[5,3,3,4,3]$ but 
we do not consider this honeycomb in this paper because hyperball covering can not be derived from them.
\end{rmrk}
\footnotesize
\medbreak
\centerline{\vbox{
\halign{\strut\vrule~\hfil $#$ \hfil~\vrule
&\quad \hfil $#$ \hfil~\vrule
&\quad \hfil $#$ \hfil\quad\vrule
&\quad \hfil $#$ \hfil\quad\vrule
&\quad \hfil $#$ \hfil\quad\vrule
\cr
\noalign{\hrule}
\multispan5{\strut\vrule\hfill\bf Table 7, 5-dimensional cases \hfill\vrule}%
\cr
\noalign{\hrule}
\noalign{\vskip2pt}
\noalign{\hrule}
 \mathcal{T} & h & Vol_5(\mathcal{S}_{pqrst}) & Vol_5(\mathcal{H}_{5}^{h}) & \delta_{pqrst}^{min} \cr
\noalign{\hrule}
[5,3,3,3,3] & 0.85377329 & 0.00076730 & 0.00133580 & 1.74091729 \cr
\noalign{\hrule}
[5,3,3,3,4] & 1.59191259 & 0.00198470 & 0.01161836 & 5.85397509 \cr
\noalign{\hrule}
\noalign{\hrule}}}}
\medbreak
\end{enumerate}
\normalsize
The next conjecture for the least dense hyperball coverings for the 3 and 5 dimensional hyperbolic spaces can be formulated:
\begin{conjecture}
{The above described hyperball covering to Coxeter tiling $[7,3,3]$ provides the least dense hyperball covering 
in the 3-dimensional hyperbolic space $\mathbf{H}^3$.}  
\end{conjecture}
\begin{conjecture}
{The above described hyperball covering to Coxeter tiling $[5,3,3,3,3]$ provides the least dense hyperball covering 
in the 5-dimensional hyperbolic space $\mathbf{H}^5$.}  
\end{conjecture}
The way of putting any analogue questions for determining the optimal ball, horoball and
hyperball
packings of tilings in hyperbolic $n$-space $(n>2)$ seems to be timely.
Our projective method suites to study and to solve these problems.

{\bf{Acknowledgement:}}
I thank Prof. Emil Moln\'ar for helpful comments to this paper.

\noindent
\footnotesize{Budapest University of Technology and Economics Institute of Mathematics, \\
Department of Geometry, \\
H-1521 Budapest, Hungary. \\
E-mail:~szirmai@math.bme.hu \\
http://www.math.bme.hu/ $^\sim$szirmai}

\end{document}